\documentclass[12pt]{amsart}
 \usepackage{amsfonts,amssymb,eucal}

\usepackage{amsthm}
 \usepackage{amsmath}
\usepackage{amscd}
 \usepackage{latexsym} 
 \newcommand{\e}{{\mbox{\rm e}}} 
  \newcommand{\mc}[1]{{\mathcal{#1}}}
 \newcommand{\got}[1]{{\mathfrak{#1}}}
\newcommand{\db}[1]{{\mathbb{#1}}}

\newcommand{\pa}{\partial}

\newcommand{\R}{\ensuremath{\mathbb{R}}}
\newcommand{\C}{\ensuremath{\mathbb{C}}}

 \newcommand{\Hi}{\ensuremath{\mathcal{H}}}
 \newcommand{\Z}{\ensuremath{\mathbb{Z}}}

 \newtheorem{Remark}{Remark}

\newtheorem{Proposition}{Proposition}

 \newtheorem{lemma}{Lemma}

\def\i{\mathrm{i}}
\begin{document}

\title{The Jacobi group and the  squeezed states - some comments}

\author{Stefan  Berceanu}
\address[Stefan  Berceanu]{National
 Institute for Physics and Nuclear Engineering\\
         Department of Theoretical Physics\\
         PO BOX MG-6, Bucharest-Magurele, Romania}
\email{Berceanu@theor1.theory.nipne.ro}

\begin{abstract} The  generalized coherent states attached to the Jacobi group  realize the squeezed states.
  Imposing hermitian conjugacy  to  the generators of the Jacobi algebra, we   
find out the form of the weight function appearing in the scalar product. We show effectively
 the orthonormality of the base functions with respect to the scalar product. From 
the explicit form of the reproducing kernel, we find out the expression of the 
multiplier in a holomorphic representation of the Jacobi group.
\end{abstract}
\maketitle


\section{Introduction\label{intr}}

 In this note we continue our investigation of the properties of the Jacobi
group started in \cite{jac1,sbj} using Perelomov's   generalized coherent states \cite{perG}. 
The Jacobi group \cite{ez,bs,neeb} -- the semidirect product of the
Heisenberg-Weyl group and the symplectic group -- is an important
object in the framework of Quantum Mechanics, Geometric Quantization,
Optics~\cite{stol1,stol2,stol3,stol4,ga1,ga2,ali,kbw2}. The Jacobi group
was investigated by physicists under other various names, as
``Schr\"odinger group'' \cite{nied} or ``Weyl-symplectic group''
\cite{kbw2}. The squeezed states \cite{stol1,stol2,stol3,stol4} in
Quantum Optics represent a physical realization of the coherent states
associated to the Jacobi group.

  In \cite{jac1} we have constructed
generalized coherent states attached to the Jacobi group,
$G^J_1=H_1\rtimes \text{SU}(1,1)$, based on the homogeneous K\"ahler
manifold
$\mc{D}^J_1=H_1/\R\times\text{SU}(1,1)/\text{U}(1)=\C^1\times\mc{D}_1$.
Here $\mc{D}_1$ denotes the unit disk $\mc{D}_1=\{w\in\C||w|<1\}$, and
$H_n$ is the $(2n+1)$-dimensional real Heisenberg-Weyl (HW) group with
Lie algebra $\got{h}_n$. In \cite{jac1} we have also emphasized the
connection of our results with those of Berndt and Schmidt \cite{bs}
and K\"ahler \cite{cal}.  In \cite{sbj} we have considered coherent
states attached to the Jacobi group $G^J_n=H_n\rtimes
\text{Sp}(2n,\R)$, defined on the manifold
$\mc{D}^J_n=\C^n\times\mc{D}_n$, where $\mc{D}_n$ is the Siegel ball.

In the present note we follow the notation and convention of
\cite{jac1}. If $\pi$ is a representation of a Lie group $G$ with Lie
algebra $\got{g}$, then we denote ${\bf{X}}=d\pi (X),~ X\in\got{g}$. In
\S \ref{unuu} we recall the  the main results on holomorphic representation in differential operators of the generators of the Jacobi algebra, while \S \ref{rep1} summarizes the information on the space of functions on which the differential operators act. \S \ref{com} contains our new results summarized in the abstract.  The results contained in Proposition \ref{fur} have been announced in \cite{sac}.

\section{The Jacobi   algebra\label{unuu}}

The Jacobi algebra is defined as the the semi-direct sum of the Lie
algebra $\got{h}_1$ of the Heisenberg-Weyl Lie group and the algebra
of the group ${\text{SU}}(1,1)$, $\got{g}^J_1= \got{h}_1\rtimes
\got{su}(1,1)$.

The Heisenberg-Weyl ideal $\got{h}_1= <\i s
1+xa^{\dagger}-\bar{x}a>_{s\in\R ,x\in\C}$ is determined by the
commutation relations
\begin{equation}
\left[a,{{K}}_+\right]=a^{\dagger}~;\left[ {{K}}_-,a^{\dagger}\right]=a;
\left[ {{K}}_+,a^{\dagger}\right]=\left[ {{K}}_-,a \right]= 0 ;
 \left[ {{K}}_0  ,a^{\dagger}\right]=\frac{1}{2}a^{\dagger};
 \left[ {{K}}_0,a\right]
=-\frac{1}{2}a,
\end{equation}
where ${a}^{\dagger}$ (${ a}$) are  the boson creation
(respectively, annihilation) operators,  which verify the canonical
commutation relations $[a,a^{\dagger}]=I, ~[a,I]= [a^{\dagger},I]=0$,
and $K_{0,+,-}$ are the generators of $\text{SU}(1,1)$ which satisfy
the commutation relations:
\begin{equation}
  \left[  {K}_{0},{K}_{\pm}\right]  =\pm {K}_{\pm}\ ,\quad\left[
    {K}_{-},{K}_{+}\right]
  =2{K}_{0}\ .\label{scs1}
\end{equation}

 We take a
 representation of $G^J_1$ (cf. \cite{jac1} and 
    Proposition \ref{mm1} reproduced below)
 such that  the cyclic vector 
$e_0$  fulfills simultaneously the
conditions
\begin{equation}\label{cond}
{\bf{a}}e_0  =  0, 
~ {\bf{K}}_-e_0   =   0,~
{\bf{K}}_0e_0   =   k e_0;~ k>0,
\end{equation}
and we take $e_0=\varphi_0\otimes\phi_{k0}$. We consider for
$\text{Sp}(2,\mathbb{R})\approx \text{SU}(1,1)$ 
the  unitary irreducible positive
discrete series representation $D^+_k$ with Casimir operator $C={
  K}_{0}^{2}-{ K}_{1}^{2}-{ K}_{2}^{2}=k(k-1)$,
where $k$ is the Bargmann index for  $D^+_k$ \cite{bar47} .
The orthonormal canonical basis of the representation space of
$\text{SU}(1,1)$ consists of the
vectors
\begin{equation}
\phi_{km} = \left[  \frac{\Gamma(2k)}{m!\Gamma(2k+m)}\right]  ^{1/2}\left(
{\bf{K}}_{+}\right)  ^{m}\phi_{k0}\ ,\label{scs2} m\in\Z_+ . %
\end{equation}

Also, in the Fock space $\mc{F}$,  we have the orthonormality $ <\varphi_{n'},\varphi_n>=\delta_{nn'}$, where 
$\varphi_n=(n!)^{-\frac{1}{2}}({\bf
  a}^{\dagger})^{n}\varphi_0$. 

Perelomov coherent state vectors associated to the group Jacobi
$G^J_1$, based on the manifold $\mc{D}^J_1$, are defined as
 \begin{equation}\label{csu}
e_{z,w}=\e^{z{\bf a}^{\dagger}+w{\bf{K}}_+}e_0, ~z,w\in\C,~ |w|<1 .
\end{equation}

The general scheme \cite{sbcag,sbl,sinaia} associates to elements of the
Lie algebra $\got{g}$ first order differential operators,
 $X\in\got{g}\rightarrow\db{X}\in \mc{D}_1$, and we have \cite{jac1}

\begin{lemma}\label{mixt}The differential action of the generators
 of
the Jacobi algebra  is given by the formulas:
\begin{subequations}\label{summa}
\begin{eqnarray}
& & \bf{a}=\frac{\pa}{\pa z};~\bf{a}^+=z+w\frac{\pa}{\pa z} ;\\
 & & \db{K}_-=\frac{\pa}{\pa w};~\db{K}_0=k+\frac{1}{2}z\frac{\pa}{\pa z}+
w\frac{\pa}{\pa w};\\
& & \db{K}_+=\frac{1}{2}z^2+2kw +zw\frac{\pa}{\pa z}+w^2\frac{\pa}{\pa
w} ,
\end{eqnarray}
\end{subequations}
where $z\in\C$, $|w|<1$.
\end{lemma}

\section{A multiplier representation of the Jacobi group}\label{rep1}

We introduce the {\it displacement operator}
\begin{equation}\label{deplasare}
D(\alpha )=\exp (\alpha {\bf a}^{\dagger}-\bar{\alpha}{\bf a)}=
\exp(-\frac{1}{2}|\alpha
|^2)  \exp (\alpha {\bf a}^{\dagger})\exp(-\bar{\alpha}{\bf a}),
\end{equation}
and  the  {\it  unitary squeezed operator} of the
  $D^k_+$ representation of the group $\text{SU}(1,1)$,
 $\underline{S}(z)=S(w)$ ($w  =
\frac{z}{|z|}\tanh \,(|z|), \eta=\ln(1-|w|^2)$):
\begin{subequations}
\begin{eqnarray}
\underline{S}(z) & = &\exp (z{\bf{K}}_+-\bar{z}{\bf{K}}_-), ~z\in\C
;\label{u1} \\
S(w) & = &  \exp (w{\bf{K}}_+)\exp (\eta
{\bf{K}}_0)\exp(-\bar{w}{\bf{K}}_-), ~|w|<1 .\label{u2}
\end{eqnarray}
\end{subequations}
 We introduce also \cite{jac1}
the {\it generalized} squeezed coherent state vector
\begin{equation}\label{gen1}
\Psi_{\alpha, w}=D(\alpha )S(w) e_0,~
e_0=\varphi_0\otimes\phi_{k0}.
\end{equation}

We introduce the auxiliary operators
\begin{equation}\label{ssq}
{\bf{K}}_+  = \frac{1}{2}({\bf a}^{\dagger})^2+{\bf{K}}'_+ ,~
{\bf{K}}_-  = \frac{1}{2}{\bf a}^2+{\bf{K}}'_- ,~
{\bf{K}}_0  = \frac{1}{2}({\bf a}^{\dagger}{\bf a}+\frac{1}{2})+{\bf{K}}'_0 ,
\end{equation}
which have the properties
\begin{equation}\label{keg}
{\bf{K}}'_-e_0   =  0,~
{\bf{K}}'_0e_0 =  k'e_0; ~k=k'+\frac{1}{4}, 
\end{equation}
\begin{equation}
\left[{\bf{K}}'_{\sigma},{\bf a}\right]=\left[
{\bf{K}}'_{\sigma},{\bf a}^{\dagger}\right]=0,~\sigma =\pm ,0 ;~
 \left[ {\bf{K}}'_0,{\bf{K}}'_{\pm}\right]=\pm {\bf{K}}'_{\pm};~
 \left[ {\bf{K}}'_-,{\bf{K}}'_+\right]=2{\bf{K}}'_0 .
\end{equation}

We recall  some properties of the coherent states associated to the
group $G^J_1$  \cite{jac1}:
  
\begin{Proposition}\label{mm1}
  The generalized squeezed coherent state
 vector {\em{(\ref{gen1})}}
and Perelomov coherent 
  state vector  {\em{(\ref{csu})}}
are related by the relation
\begin{equation}\label{csv}
  \Psi_{\alpha, w}= (1-w\bar{w})^k
  \exp (-\frac{\bar{\alpha}}{2}z)e_{z,w},\text{where}~
  z=\alpha-w\bar{\alpha}.
\end{equation}

Perelomov coherent state vector {\em{(\ref{csu})}} was calculated in
{\em{\cite{jac1}}} as
\begin{equation}\label{cemaie}
e_{z,w}=E(z,w)\varphi_0 \otimes\e^{w\bf{K}'_+}\phi_{k0},
\end{equation}
\begin{equation}\label{ezw}
  E(z,w)\varphi_0=\e^{z{\bf a}^{\dagger}+\frac{w}{2}({\bf a}^{\dagger})^2}\varphi_0
  =\sum_{n=0}^{\infty}\frac{P_n(z,w)}{(n!)^{1/2}}\varphi_n,
\end{equation}

\begin{equation}\label{marea}
  P_n(z,w)=n!\sum _{p=0}^{[\frac{n}{2}]}
  (\frac{w}{2})^p\frac{z^{n-2p}}{p!(n-2p)!} .
\end{equation}
The base of functions
$f_{nks}(z,w)=<e_{\bar{z},\bar{w}},\varphi_{n}\otimes\phi_{ks}>$,
where $k=k'+1/4,~2k'= $ integer, $n,s=0,1,\cdots $, $z, w\in\C$, $|w|<1$,  consists of functions:

\begin{equation}\label{fkn}
f_{nks}(z,w)=\sqrt{\frac{\Gamma{(s+2k-1/2)}}{s!\Gamma (2k-1/2)}}w^s\frac{P_n(z,w)}{\sqrt{n!}}.
\end{equation}

The composition law in the Jacobi group $G^J_1=HW\rtimes SU(1,1)$ is
\begin{equation}\label{compositie}
(g_1,\alpha_1,t_1)\circ (g_2,\alpha_2, t_2)= (g_1\circ g_2,
g_2^{-1}\cdot\alpha_1+\alpha_2, t_1+ t_2 +\Im
(g^{-1}_2\cdot\alpha_1\bar{\alpha}_2)),
\end{equation}
where $g^{-1}\cdot\alpha =\bar{a}\alpha -b\bar{\alpha}$, and 
$g\in\text{SU}(1,1)$ is parametrized as 
\begin{equation}\label{ggg}
g=\left(\begin{array}{cc}a & b \\\bar{b} & \bar{a}\end{array}\right),
|a|^2-|b|^2=1. 
\end{equation}

Let $h = (g,\alpha )\in G^J_1$, $\pi (h)_k = T(g)_kD(\alpha )$, $g\in
SU(1,1)$, $\alpha\in \C$, and let $x=(z,w)\in
{\mc{D}^J_1}=\C\times\mc{D}_1$.  Then we have the formulas:
\begin{equation}\label{rep}
\pi(h)_k\cdot   e_{z,w}= (\bar{a}+\bar{b}w)^{-2k}\exp(-\lambda_1)e_{z_1,w_1},
\end{equation}
\begin{equation}\label{x6}
2\lambda_1=\frac{\bar{b}}{\kappa}\gamma^2+\bar{\alpha}(z+\gamma),
\end{equation}
\begin{equation}\label{xxx}
(z_1,w_1)=(\frac{\gamma}{\kappa},\frac{aw+b}{\kappa}); ~\kappa = \bar{a}+w\bar{b}; \gamma= z+\alpha-\bar{\alpha}w .
 \end{equation}
The  space of functions $\Hi_K$ attached to the reproducing kernel
$K(z,w;\bar{z}',\bar{w}'):=(e_{\bar{z},\bar{w}},e_{\bar{z}',\bar{w}'}):
\mc{D}^J_1\times \bar{\mc{D}}^J_1\rightarrow \C $:  
\begin{equation}\label{KHK}
K(z,w;\bar{z}',\bar{w}')
=(1-{w}\bar{w}')^{-2k}\exp{\frac{2\bar{z}'{z}+z^2\bar{w}'
+\bar{z}'^2w}{2(1-{w}\bar{w}')}}, 
\end{equation}
consists of square integrable functions  with respect to the scalar
product 
\begin{equation}\label{ll2}
(f,g)_k= \Lambda\int_{z,w\in\C, |w|<1}\bar{f}(z,w)g(z,w)\rho (z,w)
\mathbf{d}^2 z\mathbf{d}^2 w, ~\Lambda =  \frac{4k-3}{2\pi^2},
\end{equation}
where
\begin{equation}\label{pp1}
\rho (z,w)=\rho_0(w) F(z,w),\end{equation}
\begin{equation}\label{pp3}
 \rho_0(w)=(1-w\bar{w})^p,~ p=2k-3,
\end{equation}
\begin{equation}\label{pp2}
 F(z,w)=
\exp\left[-\frac{2|z|^2+z^2\bar{w}+\bar{z}^2w}{2(1-|w|^2)}\right].
\end{equation}

\end{Proposition}

\section{Comments}\label{com}

\begin{Remark}\label{ora1}
If the generators of the Jacobi group $G^J_1$ have the differential realization  given in (\ref{summa}), then the operator $\bf{a}$  ($\db{K}_-$) is the hermitian conjugate of the operator  $\bf{a}^{\dagger}$,  (respectively,  $\db{K}_+$), while the operator $\db{K}_0$  is self adjoint   with respect to the scalar product (\ref{ll2}). 
\end{Remark}
\newpage

{\it Proof. } 

Using integration by parts in  the equations $({\bf{a}}f,g)=(f,{\bf{a}}^{\dagger}g)$, $(\db{K}_-f,g)=(f,\db{K}_+g)$, $(\db{K}_0f,g)=(f,\db{K}_0g)$ with respect
to  the scalar product (\ref{ll2}), we get  respectively  the equations
 \begin{subequations}\label{condd}
\begin{eqnarray}
\label{cond1}
(w\frac{\pa }{\pa z} - \frac{\pa }{\pa \bar{z}})\rho & = & z \rho,\\ 
\label{cond2}
\left( w^2\frac{\pa}{\pa w}-\frac{\pa}{\pa{\bar{w}}}
\right) \rho & = & \left(\frac{z^2}{2}+pw-zw\frac{\pa }{\pa z}\right)\rho,\\
\label{cond3}2\left( w\frac{\pa }{\pa w}-\bar{w}\frac{\pa }{\pa \bar{w}}\right)\rho & = &
\left( \bar{z}\frac{\pa}{\pa\bar{z}}-z\frac{\pa}{\pa z}\right)\rho .
\end{eqnarray}
\end{subequations}
 But 
\begin{subequations}\label{condd1}
\begin{eqnarray}\label{mm11}\frac{\pa \rho}{\pa  w} &  = & -\frac{2p\bar{w}(1-w\bar{w})+\bar{z}^2+2\bar{w}|z|^2+
\bar{w}^2z^2}{2(1-w\bar{w})^2 }\rho,\\
\label{mm2}\frac{\pa \rho}{\pa z} &  = &  -\frac{\bar{z}+z\bar{w}}{1-|w|^2}\rho,
\end{eqnarray}
\end{subequations}
and it is check out  that the function  (\ref{pp1}) - (\ref{pp2}) verifies the conditions (\ref{condd}).
\begin{Remark}\label{ora2} If $w=0$ in (\ref{cond1}), then we get the solution $\rho(z) = ct \e^{-|z|^2}$ of the reproducing kernel for the coherent states associated o the Heisenberg-Weyl group, while if $z=0$ in (\ref{cond2}), (\ref{cond3}), we get the solution  $\rho (w) = ct (1-|w|^2)^p$ of  the reproducing kernel for $SU(1,1)$.  
\end{Remark}
Starting from the Segal-Bargmann-Fock realization  of the boson operators ${\bf a}\rightarrow \frac{\pa}{\pa z}$, ${\bf a}^{\dagger}\rightarrow z$,  Bargmann \cite{bar61} has determined the reproducing kernel $\rho(z,\bar{z}) = ct \e^{-|z|^2}$ from the relation  $(zf,g)=(f,\frac{\pa g}{\pa z})$. Now we apply his technique to the Jacobi group $G^J_1$ and we obtain:  
\begin{Proposition}\label{ora3}If the generators of the Jacobi group $G^J_1$ have the differential realization  given in (\ref{summa}), then the conditions $({\bf{a}}f,g)=(f,{\bf{a}}^{\dagger}g)$, $(\db{K}_-f,g)=(f,\db{K}_+g)$, $(\db{K}_0f,g)=(f,\db{K}_0g)$ with respect
to  the scalar product (\ref{ll2}) impose to the function $\rho$ to verify the equations (\ref{condd}), which admit the solution  (\ref{pp1})-(\ref{pp2}).
\end{Proposition}
{\it Proof.} We consider for the functions $\frac{\pa \rho}{\pa z}$,  $\frac{\pa \rho}{\pa \bar{z}}$ the linear  system of    equations  consisting of    equation (\ref{cond1}) and his complex conjugate. It has as solution the equation (\ref{mm2}). Now we consider the linear  system of  equations (\ref{cond2}), (\ref{cond3}) in $\frac{\pa \rho}{\pa w}$ and $\frac{\pa \rho}{\pa {\bar{w}}}$. We introduce the solution (\ref{mm2})  for    $\frac{\pa \rho}{\pa z}$, and  we obtain the equation (\ref{mm11}). But equation (\ref{mm2}) admits the solution 
\begin{equation}\label{por}\rho(z,w,\bar{z},\bar{w})= \rho_0(w,\bar{z},\bar{w})F(z,w),
\end{equation} with $F(z,w) $ given by (\ref{pp2}). Introducing (\ref{por}) in (\ref{mm11}), we get the differential equation
$$\frac{\pa \rho_0}{\pa w}=-\frac{p\bar{w}}{1-w\bar{w}}\rho_0,$$
which has the solution (\ref{pp3}). 

We shall verify explicitly that 
\begin{Proposition}\label{fur}The system of vectors (\ref{fkn}) is orthonormal with respect to the scalar product (\ref{ll2})-(\ref{pp2}), i.e.
\begin{equation}\label{bigor}
(f_{nks},f_{mkr})=\delta_{nm}\delta_{sr}, ~k>1/2.
\end{equation} 
Let now $k\in\R$, $|k|<1/2$. Let $f=\sum_{nr}a_{nr}P_nw^r$, $g=\sum_{ms}b_{ms}P_mw^s $. Then the scalar product (\ref{ll2})-(\ref{pp2}) is replaced by
\begin{equation}\label{newsc}
(f,g)=\sum_{n,r=0}\bar{a}_{nr}b_{nr}\frac{r!\Gamma(2k-1/2)}{\Gamma(r+2k-1/2)}.
\end{equation}

The reproducing kernel (\ref{KHK}) admits the series expansion 
\begin{equation}\label{ker3}
K(z,w;\bar{z},\bar{w}')
 = \sum_{n,m}f_{nkm}(z,w)\bar{f}_{nkm}(z',w'). 
\end{equation}
\end{Proposition} 
{\it Proof.} Firstly we proof the second assertion. We can write down \cite{jac1}
\begin{equation}\label{x6x}P_n(z,w) = 
(\frac{\i}{\sqrt{2}})^n
w^{\frac{n}{2}}H_n(\frac{-\i z}{\sqrt{2w}}) .\end{equation} 
We use the summation formula (see, e.g. eq. (1.110) in \cite{grad})
\begin{equation}\label{rajik}
(1-x)^{-q}=\sum_{m=0}^{\infty}\frac{x^m}{m!}\frac{\Gamma
(q+m)}{\Gamma (q)} .
\end{equation}
We use the summation relation of the Hermite polynomials
 (Mehler formula, cf. equation 10.13.22 in  \cite{bate})
\begin{equation}\label{bat}
\sum_{n=o}^{\infty}\frac{(\frac{s}{2})^n}{n!}H_n(x)H_n(y)=
\frac{1}{\sqrt{1-s^2}}\exp{\frac{2xys-(x^2+y^2)s^2}{1-s^2}},~ |s|<1 , 
\end{equation}
and we get the second part of the assertion of the Proposition. 

Let us introduce the generating function
$$G_t(z,w):=\exp(zt +\frac{1}{2}wt^2)=\sum_{p,q\ge0}\frac{z^pw^q}{2^qp!q!}t^{p+2q}=\sum_{n\ge0}
t^n\sum_{q=0}^{[n/2]}\frac{z^pw^q}{2^qq!(n-2q)!},$$
$$G_t(z,w)= \sum_{n\ge 0}\frac{t^n}{n!}P_n(z,w).$$
Now we calculate the scalar product (\ref{ll2})-(\ref{pp2}) of two generating functions $G_t$
$$(G_t,G_t)= \Lambda \int_{z,w\in\C,|w|<1}\rho_0(w)\exp[A(z,w)]\mathbf{d}^2z\mathbf{d}^2w, 
$$
where
$$A(z,w)=-\frac{2|z|^2+z^2\bar{w}+\bar{z}^2w}{2(1-|w|^2)}+zt+\bar{z}\bar{t}+
\frac{1}{2}wt^2+\frac{1}{2}\bar{w}\bar{t}^2.$$
With the  change of  variables
$$y=\frac{z+tw-\bar{t}}{\sqrt{1-|w|^2}}; ~ (z=y\sqrt{1-|w|^2}-tw+\bar{t}),$$
we get successively  $$A(z,w)=-|y|^2-\frac{1}{2}(y^2\bar{w}+\bar{y}^2w)+|t|^2,$$
$$(G_t,G_t)=\Lambda\exp(|t|^2)\int_{|w|<1}\rho_0(w)I(w)\mathbf{d}^2w,$$
$$I(w)=\int_{\C}\exp [-|y|^2-\frac{1}{2}(y^2\bar{w}+\bar{y}^2w)]\mathbf{d}^2y=\pi(1-|w|^2)^{-1/2},$$
$$(G_t,G_t)=\pi\Lambda\exp(|t|^2)\int_{w\in\C, |w|<1}(1-|w|^2)\rho_0(w)\mathbf{d}^2w .$$
For any $k'>
1/2$ we denote $G_{ts}(z,w)=G_t(z,w)w^s$, and we have
$$(G_{ts},G_{tr})=\pi \Lambda\exp(|t|^2)\int_{w\in\C, |w|<1}\bar{w}^sw^r(1-|w|^2)^{2k-3+1/2}\mathbf{d}^2w. $$
Now we change the variable: $\Re(w)= \rho \cos\theta,~\Im(w) = \rho \sin\theta$ and then we put $\rho^2=x$. We apply formula 
$$\int_0^1t^{x-1}(1-t)^{y-1}dt=B(x,y)=\frac{\Gamma (x)\Gamma(y)}{\Gamma(x+y)}, ~\Re x, \Im x >0, $$
where $x=r+1$, $y=2k-3/2$. 
Finally, we have $$(P_nw^s,P_mw^r)=\frac{r!\Gamma (2k-1/2)}{\Gamma (r+2k-1/2)}\delta_{nm}\delta_{rs},$$
i.e. equation (\ref{bigor}). 
\begin{Remark}If the reproducing kernel $K$ is given by (\ref{KHK}), then the multiplier in the representation (\ref{rep}) of the Jacobi group $G^J_1$ has the expression (\ref{rep})-(\ref{xxx}).
\end{Remark}
{\it Proof.}


We recall the relations (cf. Prop. IV.1.9 p. 104 in \cite{neeb})
\begin{equation}\label{cedr}
K(h\cdot x,h\cdot x')=\bar{J}(h,x)K(x,x')J(h,x'),
\end{equation}
where
$$\pi(g)_kf(x)=J(g^{-1},x)^{-1}f(g^{-1}\cdot x),$$
and
\begin{equation}\label{repP}
\pi(h)_k\cdot   e_x=\lambda(h,x)e_{x_1},~ x=(z,w);~ x_1=(z_1,w_1).
\end{equation}
In the case of the Jacobi group, we get
$$(e_{z_1,w_1},e_{z'_1,w'_1})=(1-\bar{w}_1w'_1)^{-2k}\exp E',$$
where
$$E'=
\frac{2\frac{\gamma'}{\kappa'}\frac{\bar{\gamma}}{\bar{\kappa}}
+\frac{aw'+b}{\kappa'}\frac{\bar{\gamma}^2}{\bar{\kappa}^2}+
\frac{\bar{a}\bar{w}+\bar{b}}{\bar{\kappa}}\frac{\gamma'^2}{\kappa'^2}}{2(1-\bar{w}_1w'_1)};~1-\bar{w}_1w'_1=\frac{1-\bar{w}w'}{\bar{\kappa}\kappa'}.$$
We do the splitting
$$E'=E+P+Q,$$
and, if we take $$P=\frac{1}{2}[\frac{\bar{b}}{\kappa'}\gamma'^2+\frac{b}{\bar{\kappa}}\bar{\gamma}^2+\bar{\alpha}(z'+\gamma')+\alpha(\bar{z}+\bar{\gamma})],$$
then, we get $Q=0$. We have
\begin{equation}\label{dubla}K(h\cdot x,h\cdot x')=\kappa'^{2k}\exp\frac{1}{2}[\frac{\bar{b}}{\kappa'}\gamma'^2+
\bar{\alpha}(z'+\gamma')]K(x,x')\exp\frac{1}{2}[\frac{b}{\bar{\kappa}}\bar{\gamma}^2+\alpha(\bar{z}+\bar{\gamma})]{\bar{\kappa}}^{2k}.
\end{equation}

So for $h=(g,\alpha), g\in\text{SU}(1,1),~x=(w,z)$,  from (\ref{dubla}),  we get

$$J(h,x)=(\bar{b}w+\bar{a})^{2k}\exp\frac{1}{2}[\frac{\bar{b}}{\kappa}\gamma^2+
\bar{\alpha}(z+\gamma)].$$

\vspace{3ex}

{\it{Acknowledgments}.}
  The author is indebted to the Organizing Committee of the  XVIII Workshop on Geometric Methods in Physics, Bia\l owie\.{z}a, Poland 2009 for the opportunity to report results at the meeting.  The author was partially supported  by the CNCSIS  Grant ``Idei'' No.  454/2009, cod ID-44.



\begin{thebibliography}{99}







\bibitem{jac1}
S.~Berceanu, \emph{Rev. Math. Phys.} \textbf{18}, 163--199 (2006),
  {arXiv: math.DG/0408219}.

\bibitem{sbj}
S.~Berceanu, {\it{A holomorphic representation of Jacobi algebra in several
  dimensions,}} in \emph{Perspectives in Operator Algebra and Mathematical
  Physics}, edited by F.~Boca, R.~Purice, and S.~Stratila, The Theta
  Foundation, Bucharest, 2008, pp. 1--25, {arXiv: math.DG/060404381}.

\bibitem{perG}
A.~M. Perelomov, \emph{Generalized Coherent States and their Applications},
  Springer, Berlin, 1986.

\bibitem{ez}
M.~Eichler, and D.~Zagier, {\it{The theory of Jacobi forms,}} in
  \emph{Progress in Mathematics}, Birkh\"auser, Boston, MA, 1985, vol.~55.


\bibitem{bs}
R.~Berndt, and R.~Schmidt, {\it{Elements of the representation theory of
  the Jacobi group,}} in \emph{Progress in Mathematics}, Birkh\"auser Verlag,
  Basel, 1998, vol. 163.

\bibitem{neeb}K.-H.~Neeb, {\it{Holomorphy and Convexity in Lie Theory}},
 Walter de Gruyter, 2000.


\bibitem{stol1}
P.~Stoler, \emph{Phys. Rev. D} \textbf{1}, 3217--3219 (1970).

\bibitem{stol2}
E.~Y.~C. Lu, \emph{Lett. Nuovo. Cimento} \textbf{2}, 1241--1244 (1971).

\bibitem{stol3}
H.~P. Yuen, \emph{Phys. Rev. A} \textbf{13}, 2226--2243 (1976).

\bibitem{stol4}
J.~N. Hollenhors, \emph{Phys. Rev. D} \textbf{19}, 1669--1679 (1979).

\bibitem{ga1}
V.~Guillemin, and S.~Sternberg, \emph{Geometric Asymptotics}, American
  Mathematical Society, Providence, R. I., 1977.

\bibitem{ga2}
V.~Guillemin, and S.~Sternberg, \emph{Symplectic Techniques in Physics},
  Cambridge University Press, Cambridge, 1984.

\bibitem{ali}
S.~T. Ali, J.~P. Antoine, and J.-P. Gazeau, \emph{Coherent states, wavelets,
  and their generalizations}, Springer-Verlag, New York, 2000.

\bibitem{kbw2}
K.~B. Wolf, \emph{Geometric Optics on Phase Space}, Springer, 2004.

\bibitem{nied}
U.~Niederer, \emph{Helv. Phys. Acta.} \textbf{45}, 802--810 (1972/73).



\bibitem{cal}
R.~Berndt, and O.~Riemenschneider, editors, \emph{Erich K\"ahler: Mathematische
  Werke, Mathematical Works}, Walter de Gruyter, Berlin-New York, 2003.


\bibitem{sac}
S.~Berceanu, and A.~Gheorghe, \emph{Romanian J. Phys.} \textbf{53}, 1013--1021 (2008), {arXiv:  0812.0448 (math.DG)}


\bibitem{bar47}
V.~Bargmann, \emph{Ann. Math.} \textbf{48}, 568--640 (1947).



\bibitem{sbcag}
S.~Berceanu, and A.~Gheorghe, \emph{J. Math. Phys.} \textbf{33}, 998--1007
  (1992).

\bibitem{sbl}
S.~Berceanu, and L.~B. de~Monvel, \emph{J. Math. Phys.} \textbf{34}, 2353--2371
  (1993).

\bibitem{sinaia}
S.~Berceanu, {\it{Realization of coherent state algebras by differential
  operators,}} in \emph{Advances in Operator Algebras and Mathematical Physics},
  edited by F.~Boca, O.~Bratteli, R.~Longo, and H.~Siedentop, The Theta
  Foundation, Bucharest, 2005, pp. 1--24, {arXiv: math.DG/0504053}.


\bibitem{bar61}
V.~Bargmann,  \emph{Commun. Pure Appl. Math.} \textbf{24}, 187--214  (1961).


\bibitem{grad} I.~S.~Grad\v stein i  I. M. Ry\v zik, 
\emph{Tables of integrals,
sums, series and products}, (Russian),  Gosudarstv. Izdat. Fiz.-Mat. Lit., 
Moscow (1963). 
\bibitem{bate}
H.~Bateman, \emph{Higher transcendental functions}, Volume 2,
Mc Graw-Hill book,  New  York (1958).


\end{thebibliography}
\end{document}